\documentclass[11pt]{amsart}

\usepackage{amssymb,amscd,amsthm,amsxtra}
\usepackage{latexsym}

\usepackage{color}
\definecolor{citation}{rgb}{0.2,0.6,0.2}
\definecolor{formula}{rgb}{0.1,0.2,0.6} 
\definecolor{url}{rgb}{0,0,0.4}
\usepackage[colorlinks=true,linkcolor=formula,urlcolor=url,citecolor=citation]{hyperref}

\newlength{\defbaselineskip}
\newcommand{\setlinespacing}[1]
        {\setlength{\baselineskip}{#1 \defbaselineskip}}

\theoremstyle{plain}
\newtheorem{theorem}{Theorem}[section]
\newtheorem{lemma}[theorem]{Lemma}
\newtheorem{proposition}[theorem]{Proposition}
\newtheorem{corollary}[theorem]{Corollary}

\newtheorem{definition}[theorem]{Definition}
\theoremstyle{definition}
\theoremstyle{remark}
\numberwithin{equation}{section}

\newcommand{\R}{{\mathbb R}}
\newcommand{\N}{{\mathbb N}}

\newcommand{\ue}{u_{\varepsilon}}

\newcommand{\bu}{\bar{u}_{\varepsilon}}
\newcommand{\mbu}{\bar{\mu}_{\varepsilon}}

\newcommand{\tows}{\stackrel{\ast}{\rightharpoonup}}

\newcommand{\tow}{\rightharpoonup}
\newcommand{\towt}{\stackrel{\tau}{\rightarrow}}

\newcommand{\Om}{\Omega}
\newcommand{\Omb}{\overline{\Omega}}

\def\dys{\displaystyle}

\def\eps{\varepsilon}

\def\mea{\mathcal{M}({\mathbb{R}^N})}

\def\ue{u_{\eps}}
\def\mue{\mu_{\eps}}

\def\gamp{\Gamma^{+}}

\title
[Subcritical approximation of a Yamabe type non local equation]
{Subcritical approximation of a Yamabe type non local equation: a Gamma-convergence approach} 

\author[G. Palatucci]
{Giampiero Palatucci}
\address[Giampiero Palatucci]{MIPA, Universit\'e de N\^imes
\\ Site des Carmes - Place Gabriel P\'eri 30021 N\^imes, France  \quad 
\and \ Dipartimento di Matematica e Informatica, Universit\`a degli Studi di Parma
\\ Campus - Parco Area delle Scienze,~53/A
\\ 43124 Parma, Italia}
\email{\href{mailto:giampiero.palatucci@unimes.fr}{giampiero.palatucci@unimes.fr}}

\author[A. Pisante]
{Adriano Pisante}
\address[Adriano Pisante]{Dipartimento di Matematica, Sapienza Universit\`a di Roma
\\ P. le Aldo Moro, 5 
\\ 00185 Roma, Italia}
\email{\href{mailto:pisante@mat.uniroma1.it}{pisante@mat.uniroma1.it}}

\author[Y. Sire]
{Yannick Sire}
\address[Yannick Sire]{LATP %{Laboratoire d'Analyse Topologie Probabilit\'es
\\ Universit\'e Aix-Marseille III
\\ FST St. J\'er\^{o}me
\\ 13013 Marseille, France}
\email{\href{mailto:sire@cmi.univ-mrs.fr}{sire@cmi.univ-mrs.fr}}

\begin{document}

%% Ann. Scuola Norm. Sup. Pisa Cl. Sci. (5)

\begin{abstract}
We investigate a natural approximation by subcritical Sobolev embeddings of the Sobolev quotient for  the fractional Sobolev spaces~$H^s$ for any $0<s<N/2$, using $\Gamma$-convergence techniques. We show that, for such approximations, optimal functions always exist and exhibit a concentration effect of the~$H^s$ energy at one point.
\end{abstract}

\subjclass[2010]{Primary  35J60, 35C20, 35B33, 49J45\vspace{1mm}}

\keywords{Concentration-compactness principle, critical Sobolev exponent, fractional Sobolev spaces, nonlocal variational problems, $\Gamma$-convergence\vspace{1mm}}

\thanks{{\it Thanks.}\ The first author has been supported by the \href{http://prmat.math.unipr.it/~rivista/eventi/2010/ERC-VP/}{ERC grant 207573 ``Vectorial Problems''}}

\begin{center}
\begin{large}
{\sc To appear in \ Ann. Sc. Norm. Super. Pisa Cl. Sci. (5)}
\end{large}
\end{center}
\vspace{1cm}

\maketitle

\setlinespacing{1.01}

%%%%%%%%%%%%%%%%%%%%%%%%%%%%%%%%%%%%%%%%%%%%%%%%%%
%%%%%
%%%%%                     INTRODUZIONE
%%%%%
%%%%%%%%%%%%%%%%%%%%%%%%%%%%%%%%%%%%%%%%%%%%%%%%%%

\section{Introduction}

{F}or any real $s> 0$, consider the standard {\em fractional Sobolev space} $H^s$ defined via the Fourier transform
$$
H^s(\mathbb{R}^N)=\big\{ u\in L^2(\mathbb{R}^N) \, \, : \, \, |\xi|^{s}\hat{u}(\xi) \in L^2(\mathbb{R}^N)\, \big\},
$$
where 
$$
\dys \hat{u}(\xi) \, = \, \mathcal{F}(u)(\xi) \, =\, \frac{1}{(2\pi)^{\frac{N}2}}\int_{\R^N}e^{-ix\cdot\xi}u(x)\,dx.
$$
As usual, the space $H^s(\mathbb{R}^N)$ can be equivalently defined as the completion of $C_0^\infty(\mathbb{R}^N)$ with respect to the norm 
\begin{equation}
\label{Hs-norm}
\| u\|^2_{H^s}=\| (Id-\Delta)^{\frac{s}{2}} u\|^2_{L^2}=\int_{\mathbb{R}^N}(1+|\xi|^2)^s|\hat{u}(\xi)|^2 d\xi \, , 
\end{equation}
where the operator $(Id-\Delta)^{\frac{s}{2}}=\mathcal{F}^{-1}\circ M_{(1+|\xi|^2)^{s/2}} \circ \mathcal{F}$ is conjugate to the multiplication operator on $L^2(\mathbb{R}^N)$ given by the function $(1+|\xi|^2)^{s/2}$.
\vspace{2mm}

It is well known that for $0<s<N/2$,  the following Sobolev inequality does hold
for some positive constant~$S^{\ast}$, depending only on $N$ and $s$,
\begin{equation}\label{eq_sobolev}
\|u\|^{2^{\ast}}_{L^{2^\ast}} \leq S^{\ast}\|(-\Delta)^{\frac{s}{2}}u\|^{2^{\ast}}_{L^2} \ \ \ \forall u \in C^\infty_0(\mathbb{R}^N),
\end{equation}
where $2^{\ast}=2^\ast(N,s):=2N/(N-2s)$ is the Sobolev critical exponent; 
and the same inequality is valid by density on $H^s(\mathbb{R}^N)$.
\vspace{2mm}

In order to discuss inequality \eqref{eq_sobolev}, it is very natural to introduce 
the {\it homogeneous Sobolev space $\dot{H}^s$} by
\begin{equation*}
\dys \dot{H}^s(\mathbb{R}^N)=\big\{ u\in L^{2^{\ast}}\!(\mathbb{R}^N) \, \, : \, \, |\xi|^{s}\hat{u}(\xi) \in L^2(\mathbb{R}^N)\, \big\}.
\end{equation*}
This space can be equivalently defined as the completion of $C_0^\infty(\mathbb{R}^N)$ under the norm 
\begin{equation}
\label{def_Hs0norm}
\| u\|^2_{\dot{H}^s}=\| (-\Delta)^{\frac{s}{2}} u\|^2_{L^2}=\int_{\mathbb{R}^N}|\xi|^{2s}|\hat{u}(\xi)|^2 d\xi 
\end{equation}
and inequality \eqref{eq_sobolev} holds by density on $\dot{H}^s(\mathbb{R}^N)$.

\vspace{2mm}

When $0<s<1$,  an alternative formula for the norm on $\dot{H}^s(\mathbb{R}^N)$ can be given via the Gagliardo seminorm:
\begin{equation}
\label{gagliardo}
\int_{\mathbb{R}^N}|\xi|^{2s}|\hat{u}(\xi)|^2 d\xi=c(N,s) \int_{\mathbb{R}^N} \int_{\mathbb{R}^N} \frac{|u(x)-u(y)|^2}{|x-y|^{N+2s}}\, dx\,dy \, . 
\end{equation}
This can be proved by a direct calculation using Fourier transform (see, e.~\!g., \cite[Proposition~3.4]{DPV12}). However, for $s\geq 1$, the previous equality fails, since in such case the right hand-side in \eqref{gagliardo} is known to be finite if and only if $u$ is constant. 
Again in such a restrict range of validity for $s$, one can also consider the Sobolev inequality~\eqref{eq_sobolev} as the following trace inequality
\begin{equation}
\label{tracesobolev}
\dys 
\| u\|^2_{L^{2^*}}  
\leq
{S^{\ast}}^{2/2^*}\!\! \int_{\mathbb{R}^N}|\xi|^{2s}|\hat{u}(\xi)|^2 d\xi  \leq 
C(N,s) \int_{\mathbb{R}^N}\int_0^\infty |\nabla U|^2 t^{1-2s} \,dx\,dt \, ;
\end{equation}
see~\cite{CS07}.

\vspace{2mm}

In this paper we are interested into investigate a natural approximation to~\eqref{eq_sobolev} via subcritical Sobolev inequalities, from a variational point of view, in the full range of validity $0< s< N/2$. For this, from now on we have to deal with the norm using Fourier as in~\eqref{def_Hs0norm}. Also, it is convenient to consider   
the following maximization problem
\begin{eqnarray}\label{pbsobolev}
\dys S^{\ast} \! := \!\sup\left\{ \int_{\mathbb{R}^N} |u|^{2^{\ast}\!} dx :  u\in \dot{H}^s(\mathbb{R}^N), \, \int_{\R^N} |(-\Delta)^{\frac{s}{2}} u|^{2} dx \leq 1 \right\}.
\end{eqnarray}
Clearly, the validity of~\eqref{eq_sobolev} is equivalent to show that the constant $S^\ast$ given by~\eqref{pbsobolev} is finite. The explicit form of the maximizers, whose existence is not trivial since the critical embedding is not compact because of dilation and translation invariance, has been shown by Cotsiolis \& Tavoularis in~\cite{cotsiolis}, together with the computation of the optimal constant; see Section~\ref{sec_preliminaries}.

\vspace{2mm}

Analogously, for any bounded domain $\Omega\subset\R^N$, 
one can consider the following 
 maximization problem (or Sobolev embedding)
\begin{eqnarray}\label{pbsobolev2}
\dys S^{\ast}_\Omega \! := \!\sup\left\{\int_{\Omega} |u|^{2^{\ast}\!} dx: u\in \dot{H}^s(\Omega), \, \int_{\R^N} |(-\Delta)^{\frac{s}{2}} u|^{2} dx \leq 1 \right\},
\end{eqnarray}
where the Sobolev space $\dot{H}^s(\Omega)$ is given as the closure of $C^\infty_0(\Omega)$ in $\dot{H}^s(\mathbb{R}^N)$ with the norm in \eqref{def_Hs0norm}.

A simple scaling argument on compactly supported smooth functions shows that $S^{\ast}=S^{\ast}_\Omega$, but in view of~\cite[Theorem 1.1]{cotsiolis}, the variational problem \eqref{pbsobolev2} has no maximizer (see Theorem~\ref{the_cotsiolis} below).
This is not the case when one considers the subcritical embeddings. {F}or any $0<\eps < 2^\ast-2$, we set
\begin{eqnarray}\label{problemaeps}
\dys S^{\ast}_{\eps} \! := \!\sup\left\{ F_\eps(u) : u \in \dot{H}^s(\Omega), \, \int_{\R^N} |(-\Delta)^{\frac{s}{2}} u|^{2} dx \leq 1 \right\}.
\end{eqnarray}
\begin{eqnarray}\label{feps} 
\dys \text{where} \ \ F_\eps (u)\!\! := \!\!\int_{\Om} |u|^{2^{\ast}\!-\eps} dx.  
\end{eqnarray}
Clearly, since $\Omega$ is a bounded domain, the embedding $\dot{H}^s(\Om) \hookrightarrow L^{2^\ast-\eps}(\Om)$ is compact and this will assure the existence of a maximizer $u_\eps\in \dot{H}^s$ for the previous problem.

\vspace{2mm}

The aim of this paper is to investigate what happens 
when $\varepsilon~\to~0$ both to the subcritical Sobolev constant $S^{\ast}_\eps$ given in~\eqref{problemaeps} and to the corresponding maximizers~$u_\varepsilon$, i.~\!e. the corresponding optimal functions of the embedding $\dot{H}^s(\Omega) \hookrightarrow L^{2^{\ast}-\varepsilon}(\Omega)$.

\vspace{2mm}

At least in the {\it local case} $s=1$, this problem has been widely investigated during the last decades, mainly, by studying the Euler-Lagrange equation for the functional $F_\varepsilon$ given by~\eqref{feps}, among functions with $\dot{H}^1$ norm equal to one; that is,  
\begin{equation}\label{eq_diri}
-\Delta u_{\eps} = \lambda |u_{\eps}|^{2^{\ast}-2-\eps}\ue \ \ \text{in} \ (\dot{H}^1(\Omega))',
\end{equation}
where $\lambda$  
is a Lagrange multiplier. 

In particular, we would mention the preliminary paper~\cite{atkinsonpeletier} in the case when $\Om$ is the unit ball in $\R^3$, where it has been showed that the solutions $u_\eps$ of \eqref{eq_diri}, maximizing the Sobolev quotient, are such that
$$
\dys \lim_{\eps\to 0} \eps u_\eps^2 (0) = \frac{32}{\pi} \ \ \ \text{and} \ \ \ \lim_{\eps\to0}\eps^{-1/2}\ue(x)=\frac{\sqrt{\pi}}{4\sqrt{2}}\left(\frac{1}{|x|}-1\right), \ \forall x\neq 0.
$$
In   \cite{brezispeletier}, Brezis \& Peletier extended such result to the case of $\Om$ being a spherical domain, along with other interesting statements. In particular, in \cite{brezispeletier} it has been  showed that the subcritical solutions concentrate at one special point of $\Om$; the authors also conjectured that the same kind of results holds for non spherical domains.
Such conjecture has been proved in the case of any smooth bounded domain $\Om$ by Han in \cite{han} and by Rey in \cite{rey}, by showing that the solutions of (\ref{eq_diri}), with maximal Sobolev energy,
concentrate at one point $x_0\in \Om$  being a critical point  for the Robin function $\mathcal{R}_\Om$, the diagonal of the regular part of the Green function in $\Om$.

\vspace{1mm} 

In order to obtain this concentration result, even without localizing the blowing-up, the proofs of all the  results above also utilize  standard elliptic regularity techniques that require to work in smooth domains.
We would notice that further regularity assumptions on $\Om$ are not for free, since they will imply that the associated Robin function $\mathcal{R}_\Om$ diverges at the boundary, and this yields that the concentration point $x_0$ belongs to the interior of $\Om$. This is not the case in~\cite{pistoia}, when a concentration up to the boundary has been proven in an example of a nonsmooth domain $\tilde{\Omega}$, previously defined in~\cite{fluchergarroni}.

\vspace{1mm}

Recently, one of the authors has proven similar concentration results  when $s=1$   in the case of any bounded domain~$\Om$ with no regularity assumptions, as well as describing the asymptotic analysis of the Sobolev quotient~\eqref{problemaeps} in term of De~Giorgi's~$\Gamma$-convergence (see~\cite{palatucci09}, and also~\cite{palatucci10} for the nonlinear case in the Sobolev spaces $W^{1,p}$).

\vspace{2mm}

Here, we can extend the analysis of the subcritical Sobolev embeddings in the more general $\dot{H}^s$ framework, for any $0<s<N/2$, in turn implying the natural concentration result in any bounded domain, possibly non smooth.  
 In this direction, our analysis is one of the first in the nonlocal framework, together with the one related to some improved fractional Sobolev embeddings and to the concentration-compactness principle for bounded sequences in the fractional Sobolev spaces, established by two of the authors in the forthcoming paper~\cite{PP13}.
Also, it is worth mentioning the relevant papers on the fractional Yamabe problem on manifold with boundary by Escobar~(\cite{Es}) in the case when $s=1/2$, and by Gonzalez and Qing~(\cite{GQ}) when~$s\in (0,1)$.

\vspace{1mm}
Let us come back to the variational form in the subcritical problem~\eqref{problemaeps}. We will carefully analyze the asymptotic behavior of the energy functionals $F_\eps$, by means of $\Gamma$-convergence techniques.  
The analysis here is much in the spirit of~\cite{palatucci09,palatucci10} and~\cite{amargarroni},
but with serious differences in the proofs, firstly due to the nonlocality of the fractional Sobolev spaces (for which we refer to~Section~\ref{sec_lemma}). Moreover,  
in the aforementioned papers the existence of a recovery sequence is proven by means of compactness and locality properties of the $\Gamma$-limit. Here, we pursue a different strategy that will permit us to explicitly exhibit   such recovery sequences;
see also the remarks at the beginning of Section~\ref{sec_gammaprof}.

\vspace{1mm}

Our main result is condensed in the following

\begin{theorem}\label{the_gamma-intro}
Let $\Om\subseteq\R^N$ be a bounded domain and let $X$ be the space 
$$
X=X(\Om):=\bigl\{(u,\mu) \in \dot{H}^s(\Omega)\times\mea: \mu \geq |(-\Delta)^{\frac{s}{2}} u|^{2}dx, \, \mu({\R^N})\leq 1\bigl\},
$$ endowed with the product  topology $\tau$ such that
\begin{equation}\label{def_top1}
(u_n,\mu_{n}) \towt (u,\mu) \ \, {\stackrel{\text{def}}{\Leftrightarrow}} \ \, \begin{cases} u_{n} \rightharpoonup u \ \text{in} \ L^{2^{\ast}}\!(\Omega), \\
\mu_{n} \tows \mu \ \text{in} \ \mea.
\end{cases}
\end{equation} 
\noindent
\\ Let us consider the following family of functionals 
\begin{equation}\label{def_fue1}
{F}_{\eps}(u,\mu):= \dys \int_{\Omega}|u|^{2^{\ast}\!-\eps} dx  \ \ \forall (u,\mu) \in X
\end{equation}
Then, as $\eps\to 0$, the $\Gamma^{+}$-limit of the family of functionals ${F}_{\eps}$ with respect to the topology $\tau$ corresponding to {\rm(\ref{def_top1})} is the functional ${F}$ defined by
$$
{F}(u,\mu)=\int_{\Omega}|u|^{2^{\ast}}dx + S^{\ast}\sum_{j=1}^{\infty}\mu_{j}^{\frac{2^{\ast}}{2}}, \ \ \ \forall (u,\mu) \in X.  
$$
Here $S^{\ast}$ is the best Sobolev constant in $\R^N$ 
and the numbers $\mu_j$ are the coefficients of the atomic part of the measure $\mu$. 
\end{theorem}
It would be interesting to prove results analogous to those in the preceding theorem with respect to the equivalent norms \eqref{gagliardo} and \eqref{tracesobolev}, 
 i.~\!e., taking the measure $\mu$ as limit of energy densities in $\mathbb{R}^N\times \mathbb{R}^N$ or $\mathbb{R}^N \times (0,\infty)$, respectively, and describing the corresponding loss of compactness in terms of atomic measures (as in Theorem~\ref{the_cca} below).

As a consequence of the result in Theorem~\ref{the_gamma-intro}, together with 
  the $\Gamma^{+}$-convergence property of
convergence of maximizers, 
we can also deduce that the sequences of maximizers~$\{u_\varepsilon\}$ for $S^{\ast}_\eps$ concentrate energy at one point $x_0\in\Omb$, in clear accordance with the local case.
\begin{corollary}\label{cor_concentration}
Let $\Omega \subset\mathbb{R}^N$ be a bounded open set and for each $0<\eps<2^\ast-2$ let $\ue\in \dot{H}^s(\Omega)$ be a maximizer for $S^{\ast}_{\eps}$.  Then, 
as $\varepsilon=\varepsilon_n \to 0$, up to subsequences $u_n={\ue}_n$ satisfies $u_n \rightharpoonup 0$ in $L^{2^\ast}\!(\Omega)$ and it concentrates at some point $x_{0}\in\Omb$ in~$H^s$, i.e.
$$
\dys  |(-\Delta)^{\frac{s}{2}}u_n|^{2}dx \tows \delta_{x_{0}} \ \text{in} \ \mathcal{M}(\mathbb{R}^N).
$$
\end{corollary}

\vspace{1mm}

Finally, it is worth noticing that Theorem~\ref{the_gamma-intro} could have its own relevance also to identify the location of the concentration point in the nonlocal case. Indeed, it can be read as the necessary first step in the asymptotic development by $\Gamma$-convergence (as firstly introduced in~\cite{AB93}; see also the recent paper~\cite{BT08}) of the functionals in~\eqref{def_fue1}. In this sense, a second order expansion of the $\Gamma$-limit could bring the desired informations on the concentration of the maximizing sequences, as in~\cite{amargarroni}, where different energies involving critical growth problems have been studied (see, also,~\cite{fluchergarroni} and~\cite{garronimuller}).

\vspace{3mm}

In this respect, another subcritical problem that would be very natural to investigate is the fractional counterpart of the Brezis-Nirenberg problem
\begin{equation}\label{brezisnirenberg}
(-\Delta)^{s} u -\eta u=  |u|^{2^{\ast}-2}u \ \ \text{in} \ (\dot{H}^s(\Omega))',
\end{equation}
where $\eta>0$ is a parameter. Well known results for $s=1$ (see \cite{BrNi}) and $s=2m$ an even integer (see~\cite{PS}) suggest that, even for fractional values of $s$, existence results for \eqref{brezisnirenberg} should always depend in a delicate way on $\eta$ (for preliminary results in this direction, see, e.~g.,~\cite{SV13} when $s\in (0,1)$, and~\cite{BCD11,tan} though with a slightly different definition of the fractional Laplacian;
see also~\cite{DPV12b} and \cite{Sec12} 
for related results for nonlinear fractional Schr\"odinger equations).

\vspace{3mm}

The paper is organized as follows. In Section~\ref{sec_preliminaries} below, we recall some recent results about the functions in the fractional Sobolev space~$\dot{H}^s$, as the concentration-compactness alternative and some workarounds to handle the nonlocality of the fractional Laplacian, which will be relevant in the rest of the paper. Section~\ref{sec_gamma} is devoted to the $\Gamma$-convergence analysis as given by Theorem~\ref{the_gamma-intro} and  subsequently to the proof of the concentration result in Corollary~\ref{cor_concentration}.

\vspace{4mm}

%%%%%%%%%%%%%%%%%%%%%%%%%%%%%%%%%%%%%%%%%%%%%%%%%%%
%
%               PRELIMINARY RESULTS
%
%%%%%%%%%%%%%%%%%%%%%%%%%%%%%%%%%%%%%%%%%%%%%%%%%%%

\section{Preliminaries}\label{sec_preliminaries}

In this section, we recall some recent results involving the fractional Sobolev inequality~\eqref{eq_sobolev} and the  analysis of the effects of the corresponding lack of compactness.

\vspace{2mm}

Firstly, we state the aforementioned theorem proved in~\cite{cotsiolis} which gives the optimal constant in the Sobolev inequality~\eqref{eq_sobolev} together with the explicit formula for those functions giving equality in the inequality.

\begin{theorem}{\rm (\cite[Theorem 1.1]{cotsiolis})}\label{the_cotsiolis}
Let $0<s<N/2$ and $2^{\ast}=2N/(N-2s)$. Then
\begin{equation}\label{eq_sobolev1}
\|u\|^{2^{\ast}}_{L^{2^{\ast}}\!(\R^N)} \leq S^{\ast}\|(-\Delta)^{\frac{s}{2}}u\|^{2^{\ast}}_{L^2(\R^N)} \ \ \ \forall u \in \dot{H}^s(\R^N),
\end{equation}
where
$$
S^{\ast}=\left(2^{-2s}\pi^{-{s}}\frac{\mathbf\Gamma\left(\frac{N-2s}{2}\right)}{\mathbf\Gamma\left(\frac{N+2s}{2}\right)}\left[\frac{\mathbf\Gamma(N)}{\mathbf\Gamma(N/2)}\right]^{2s/N}\right)^{\!\!\frac{2^{\ast}}{2}}\!\!
$$
and $\mathbf{\Gamma}$ is the Gamma function.
\\ For $u\neq 0$, we have equality in {\rm (\ref{eq_sobolev1})} if and only if
\begin{equation}\label{def_talentiana}
\dys
u(x)=\frac{c}{(\lambda^2+|x-x_0|^2)^{\frac{N-2s}{2}}} \ \ \forall x\in \R^N,
\end{equation}
where $c \in \R\setminus\{0\}$, $\lambda >0$ and  $x_0\in \R^N$ are fixed constants.
\end{theorem}
%\vspace{2mm}
%
The theorem above in the case~$s=1$ has been proved in~\cite{talenti} and also in~\cite{aubin}, where the connection with the Yamabe problem is also discussed. When $2\leq s<N/2$ is an even integer the same result was obtained some years later in~\cite{Sw}, following the ideas in~\cite{lions} and~\cite{lions2}.
The proof in \cite{cotsiolis} is based on a sharp form of the Hardy-Littlewood-Sobolev inequality.
By means of the moving planes method, formula \eqref{def_talentiana} has been also obtained by Chen, Li \& Ou in \cite{ouLi}; and, at least when $0<s<1$, in~\cite{FS} one can find a third approach through symmetrization techniques applied to the norm in the right hand-side of~\eqref{gagliardo}.

\vspace{2mm}

An important contribution in order to study the behavior of a maximizing sequence for \eqref{pbsobolev} and \eqref{pbsobolev2} is to establish a concentration-compactness alternative for bounded sequences in the fractional space $\dot{H}^s$, as stated in the following theorem, proved in~\cite{PP13} using ideas and methods introduced in the pioneering works \cite{lions} and \cite{lions2}.
\begin{theorem}\label{the_cca}
{\rm (\cite[Theorem 1.5]{PP13}).} Let $\Omega \subseteq \mathbb{R}^N$ an open subset and let $\{u_n\}$ be a sequence in $\dot{H}^s(\Om)$ weakly converging to $u$ as $n \to \infty$ and such that
$$
|(-\Delta)^{\frac{s}{2}} u_n|^2dx \tows \mu \ \ \ \text{and} \ \ \ |u_n|^{2^{\ast}}dx\tows \nu \ \ \text{in} \ \mea.
$$ 
Then, either $u_n \to u$ in $L^{2^{\ast}}_{\rm{loc}}(\mathbb{R}^N)$ or  there exists a (at most countable) set of distinct points $\{x_j\}_{j\in J}$ and positive numbers $\{\nu_j\}_{j\in J}$ such that we have
\begin{equation}
\label{quantnu}
\nu=\ |u|^{2^{\ast}}dx+\sum_{j} \nu_j \delta_{x_j}.
\end{equation}
If, in addition, $\Omega$ is bounded, then there exist a positive measure $\tilde{\mu} \in \mathcal{M}(\mathbb{R}^N)$ with {\rm spt}~$\!\tilde{\mu}~\subset~\Omb$ and positive numbers $\{\mu_j\}_{j\in J}$  such that
\begin{equation}
\label{quantmu}
\mu=|(-\Delta)^{\frac{s}{2}}u|^2dx+\tilde{\mu}+\sum_{j} \mu_j \delta_{x_j}, \quad \nu_j \leq S^{\ast} (\mu_j)^{\!\frac{2^{\ast}}2}\, .
\end{equation}
\end{theorem}
A consequence of the previous theorem  which will be useful in the next section is the following result which shows that on bounded domains there is no energy loss in the concentration process. 

\begin{proposition}
\label{pro_supporto}
{\rm (\cite[Proposition 6.5]{PP13}).}Let $0<2s<N$, let  $\Omega \subset \mathbb{R}^N$ be a bounded open set and let $\{u_n\}\subset \dot{H}^s(\Omega)$ such that $u_n \rightharpoonup 0$ as $n \to \infty$. For any open set $A\subseteq \mathbb{R}^N$ such that $\overline{\Omega} \cap \overline A=\emptyset$ we have $\dys \int_A |(-\Delta)^{\frac{s}{2}} u_n|^2 dx \to 0$ as $n \to \infty$.
\end{proposition}

\vspace{2mm}

We conclude this section with
\subsection{Some useful lemmas}\label{sec_lemma}

One of the main difficulties to handle functions in the fractional  Sobolev spaces is given by the intrinsic nonlocality. 
In the rest of the paper we will make use of the following two lemmas which are workarounds to use cut-off functions and provide a way to manipulate smooth truncations for the fractional Laplacian; their proofs carefully requires properties of multipliers between Sobolev spaces and strong commutator estimates (see~\cite{PP13}).

\begin{lemma}
\label{lem_cutoff}
{\rm (\cite[Lemma 6.1]{PP13}).}
Let $0<s<N/2$ and let  $u\in \dot{H}^s(\mathbb{R}^N)$. Let $\varphi \in C^\infty_0(\mathbb{R}^N)$ and for each $\lambda>0$ let $\varphi_\lambda (x):=\varphi(\lambda^{-1}x)$. Then
$$
u\varphi_\lambda \to 0 \ \text{in} \ \dot{H}^s(\mathbb{R}^N) \ \text{as}\ \lambda \to 0.
$$
If, in addition, $\varphi\equiv 1$ in a neighborhood of the origin, then
$$
u\varphi_\lambda \to u  \ \text{in} \ \dot{H}^s(\mathbb{R}^N) \ \text{as}  \ \lambda \to \infty.
$$
\end{lemma}

\vspace{2mm}

\begin{lemma}
\label{lem_commutator}
{\rm (\cite[Lemma 6.2]{PP13}).}
Let $0<s<N/2$, let $\Omega \subset \mathbb{R}^N$ a bounded open set and let $\varphi~\!\in~\!C^\infty_0(\mathbb{R}^N)$. Then the commutator $\, [\varphi, (-\Delta)^{\frac{s}{2}}]:\dot{H}^s(\Omega) \to L^2(\mathbb{R}^N) \,$ is a compact operator, i.e.
$$
\varphi  ((-\Delta )^{\frac{s}{2}}u_n)- (-\Delta)^{\frac{s}{2}} (\varphi u_n) \to 0 \quad \hbox{in} \quad L^2(\mathbb{R}^N)
$$
whenever $u_n \rightharpoonup 0$ in $\dot{H}^s(\Omega)$ as $n \to \infty$.
\end{lemma}

\vspace{1mm}

%%%%%%%%%%%%%%%%%%%%%%%%%%%%%%%%%%%%%%%%%%%%%%%%%%%
%
%                     GAMMA-CONVERGENCE
%
%%%%%%%%%%%%%%%%%%%%%%%%%%%%%%%%%%%%%%%%%%%%%%%%%%

\section{Subcritical approximation of the Sobolev quotient}\label{sec_gamma}

Let $\Omega$ be a bounded domain in $\mathbb{R}^{N}$, $s\in \R$,  $0<s<N/2$. In the introduction we have considered the following problem for $\varepsilon \in (0, 2^{\ast}-2)$
\begin{eqnarray}\label{problema}
\dys S^{\ast}_{\eps} \! = \!\sup\left\{\int_{\Om} |u|^{2^{\ast}\!-\eps} dx : \, u \in \dot{H}^s(\Omega), \, \int_{\R^N} |(-\Delta)^{\frac{s}{2}} u|^{2} dx \leq 1 \right\}.
\end{eqnarray}
The goal of the present section is to  describe the asymptotic behavior as $\eps \to 0$ of the optimal constants $S^\ast_\eps$ associated to the embeddings $\dot{H}^s(\Om) \hookrightarrow L^{2^{\ast}-\eps}(\Om)$ and of the corresponding maximizers, studying the family $\{F_\eps\}$ of functionals
\begin{eqnarray}\label{famiglia}
\dys F_\eps (u)\!\! := \!\!\int_{\Om} |u|^{2^{\ast}\!-\eps} dx,
\end{eqnarray}
on the set $\left\{u \in \dot{H}^s(\Omega), \, \int_{\R^N} |(-\Delta)^{\frac{s}{2}} u|^{2} dx \leq 1 \right\}$. 

The main tool is the notion of $\Gamma$-convergence in the sense of De Giorgi (see~\cite{dalmaso92} for an introduction) and the crucial point is to introduce a convenient functional framework in which performing the passage to the limit, i.~\!e., the functional space $X$ endowed with the topology $\tau$ as given in Theorem~\ref{the_gamma-intro}.

The reason for the choice of $X$ 
 can be described as follows. We are interested in the asymptotic behavior of the sequence $\{{F}_{\eps}(u_{\eps})\}$ for every sequence $\{u_{\eps}\}$ such that $\|(-\Delta)^{\frac{s}{2}} u_{\eps}\|^{2}_{L^{2}(\R^N)}\leq 1.$
The constraint on the ``Dirichlet energy'' of $u_\eps$ implies that, up to subsequences, there exists $\mu \in \mea$ and $u \in \dot{H}^s(\Om)$ such that $\mu({\R^N})\leq 1$, $|(-\Delta)^{\frac{s}{2}} u_{\eps}|^{2}dx\tows \mu$ in $\mea$ and $u_\eps \rightharpoonup u$ in $\dot{H}^s$. 
Clearly, by Sobolev embedding, we also have $u_{\eps} \rightharpoonup u$ in $L^{2^{\ast}}\!(\Om)$. \vspace{1mm}
By Fatou's Lemma, we deduce $\mu\geq|(-\Delta)^{\frac{s}{2}} u|^{2}dx$ and we can always decompose $\mu$ in
$$
\mu=|(-\Delta)^{\frac{s}{2}} u|^{2}dx+\tilde{\mu}+\sum_{j=1}^{\infty}\mu_{j}\delta_{x_{j}},
$$
where $\mu_{j} \in [0,1]$ and $\{x_{j}\} \subseteq \overline{\Omega}$ are distinct points; the positive measure $\tilde{\mu}$ can be viewed as the ``non-atomic part'' of the measure $(\mu-|(-\Delta)^{\frac{s}{2}} u|^2dx)$. \vspace{1mm}
In view of this decomposition, the definition of $X$ given in Theorem~\ref{the_gamma-intro} is very natural; moreover the space~$X$ is sequentially compact in the topology $\tau$.
Indeed, if $\{u_n, \mu_n\}\subseteq X$, then $\{u_n\}$ is bounded in $\dot{H}^s(\Om)$. Up to subsequences, $\mu_n\tows\mu$ in $\mea$ and $u_n \rightharpoonup u$ in $\dot{H}^s(\Om)$ (and in $L^{2^\ast}\!(\Om)$, by Sobolev embeddings) and the inequalities defining $X$ still hold for $(u,\mu)$ by weak lower semicontinuity.

\vspace{2mm}
Since $X$ appears as a sort of completion of $\dot{H}^s(\Om)$ in the weak topology of the product $L^{2^\ast}\!(\Om)\!\times\mea$, 
it would be interesting to understand whether, as in the case $s=1$ (see \cite[Proposition 2.3]{amargarroni}), every pair~$(u,\mu)$ in $X$ can be actually approximated in the topology~$\tau$ by a sequence of the form $\{(\ue, |(-\Delta)^{\frac{s}{2}}\ue|^2dx)\}$. We will not pursue this point here.
\vspace{1mm}

%-------------------

Note that, since the embeddings  $\dot{H}^s(\Om) \hookrightarrow L^{2^{\ast}-\eps}(\Om)$ are compact 
(see for instance~\cite[Lemma 10]{PSV12} for a simple proof),
the functionals $F_\eps$ as extended to $X$ by~\eqref{famiglia} are continuous and Proposition~\ref{pro_argmax} below show that there are no further maximizers in the space $X$.

As a consequence, we have that the $\Gamma^{+}$-convergence of functionals in this space implies the convergence of maximizers $\{\ue\}$ of $F_\eps$ to the maxima of $F$; this will allow an alternative proof of the concentration for the sequences $\{\ue\}$.

\begin{proposition}\label{pro_argmax}
For any $\eps>0$, let $(\bar{u}_\eps, \bar{\mu}_\eps)\in X$ be such that 
$$
\dys \sup_{(u,\mu)\in X} F_\eps(u,\mu) = F_\eps(\bar{u}_\eps, \bar{\mu}_\eps).
$$
Then 
$
\bar{\mu}_\eps=|(-\Delta)^{\frac{s}{2}}\bar{u}_\eps|^2dx.
$
\end{proposition}
\vspace{1mm}
\begin{proof}
We observe that the supremum is attained at some $(\bar{u}_{\eps},\mbu)$ because $X$ is sequentially compact and $F_\eps$ is sequentially continuous (due to the compact embedding $\dot{H}^s(\Om) \hookrightarrow L^{2^\ast-\eps}(\Om)$).

Clearly, we may suppose $\dys \bar{\mu}_\eps(\R^N)=1$.
Indeed, if we have $\lambda_\eps:=\bar{\mu}_\eps(\R^N)<1$, then we may consider the pair $(\bar{u}_\eps, \bar{\mu}_\eps/\lambda_\eps)$ which belongs to the space $X$ and satisfies $(\bar{\mu}_\eps/\lambda_\eps)(\R^N)=1$ and
$$
F_\eps(\bar{u}_\eps, \bar{\mu}_\eps/\lambda_\eps)=F_\eps(\bar{u}_\eps, \bar{\mu}_\eps)=\max_{(u,\mu)\in X} F_\eps(u,\mu).
$$
Since $\bar{u}_\eps\neq0$, by the definition of $X$ we have $0<\|(-\Delta)^{\frac{s}{2}}\bar{u}_\eps\|_{L^2}\leq1$.
Hence, if we set
\begin{equation}\label{eq_recall}
\alpha\, =\, \dys \alpha(\eps):=\, \frac{1}{\|\bar{u}_\eps\|^2_{\dot{H}^s(\Om)}}\, \geq \, 1,
\end{equation}
we may consider a new  pair $(\tilde{u}_\eps,\tilde{\mu}_\eps)$ given by
$$
\dys \tilde{u}_\eps:=\sqrt{\alpha}\bar{u}_\eps \ \ \text{and} \ \ \tilde{\mu}_\eps:=\alpha|(-\Delta)^{\frac{s}{2}}\bar{u}_\eps|^2dx.
$$
Note that $(\tilde{u}_\eps,\tilde{\mu}_\eps)$ belongs to the space $X$ and it satisfies
\begin{eqnarray}\label{eq_maxima}
F_\eps(\tilde{u}_\eps,\tilde{\mu}_\eps) \! & = & \! \alpha^{\frac{{2^{\ast}-\eps}}{2}}F_\eps(\bar{u}_\eps,\bar{\mu}_\eps) \, = \, \alpha^{\frac{{2^{\ast}-\eps}}{2}}\max_{(u,\mu)\in X} F_\eps(u,\mu).
\end{eqnarray}
Clearly, \eqref{eq_recall} and \eqref{eq_maxima} imply that $\alpha=1$, 
$(\bar{u}_\eps, |(-\Delta)^{\frac{s}{2}}\bar{u}_\eps|^2dx)$ is a maximizer and $\|\bar{u}_\eps\|_{\dot{H}^s(\Om)}=1$. 
Since $1=\int |(-\Delta)^{\frac{s}{2}}\bar{u}_\eps|^2dx \leq \mbu(\R^N)=1$, we have $\mbu=|(-\Delta)^{\frac{s}{2}}\bar{u}_\eps|^2dx$ and the proof is complete. 
\end{proof}

\subsection{Proof of the $\Gamma^{+}$-convergence result}\label{sec_gammaprof}
We first recall the definition of $\Gamma^{+}$-convergence adapted to our framework (see \cite{dalmaso92} for further details).
\begin{definition} 
We say that the family $\{{F}_{\eps}\} \ \Gamma^{+}$-converges to a functional ${F}: X \rightarrow [0,\infty),$ as $\eps \to 0,$ if for every $(u,\mu) \in X$ the following conditions hold:
\begin{itemize}
\item[(i)]{for every sequence $\{(u_{\eps},\mu_\eps)\} \subset X$ such that $ u_{\eps}\rightharpoonup u$ in $L^{2^{\ast}}\!(\Om)$ and $\mu_\eps\tows \mu$ in $\mea$
$$
{F}(u,\mu)\geq \limsup_{\eps \to 0}{F}_{\eps}(u_{\eps},\mu_\eps);
$$
}
\item[(ii)]{there exists a sequence $\{\bar{u}_{\eps}, \mbu)\} \subset X$ such that $\bar{u}_{\eps} \tow u$ in $L^{2^{\ast}}\!(\Om)$, $\bar \mu_\eps \tows \mu$ in $\mea$ and
$$
{F}(u,\mu)\leq \liminf_{\eps \to 0}{F}_{\eps}(\bar{u}_{\eps}, \bar \mu_\eps).
$$
}
\end{itemize}
\end{definition}
The $\Gamma^{+}$-limsup inequality~(i) easily follows from the concentration-compactness alternative  stated in Section \ref{sec_preliminaries} (see forthcoming Proposition~\ref{pro_limsup}). The proof of the $\Gamma^{+}$-liminf inequality~(ii) (i.e., the construction of a recovery sequence) is more delicate. In the case $s=1$ it is proved in~\cite{palatucci09}, following the strategy adopted in the proof of~\cite[Theorem 3.1]{amargarroni}.
As already mentioned in the introduction, 
both in~\cite{amargarroni} and in~\cite{palatucci09}, the authors prove the existence of a recovery sequence and the $\Gamma^{+}$-liminf inequality, working in two separate cases $(u,\mu)=(u,|\nabla u|^{2}dx+\tilde{\mu})$ and $(u,\mu)=\left(0, \sum_{i}\mu_{i}\delta_{x_{i}}\right)$ and cover the general case by means of compactness and locality properties of the $\Gamma^{+}$-limit. Here, we follow a different strategy and we explicitly construct a recovery sequence using the optimal functions given by Theorem~\ref{the_cotsiolis}.
\vspace{1mm}

The proof of the $\gamp$-limsup inequality (i) is given by the following result.
\begin{proposition}\label{pro_limsup}
For every $(u,\mu) \in X$ and for every sequence $\{(u_{\eps},\mu_{\eps})\} \subset X$ such that $(u_{\eps},\mu_{\eps})\towt (u,\mu)$, we have
$$
{F}(u,\mu) \geq \limsup_{\eps \to 0}{F}_{\eps}(u_{\eps},\mu_{\eps}).
$$
\end{proposition}
\vspace{1mm}
\begin{proof}
Let $\{(u_{\eps},\mu_{\eps})\}$ be a sequence in $X$ such that $(u_{\eps},\mu_{\eps})\towt(u,\mu)$;
clearly, $\dys \mu=|(-\Delta)^{\frac{s}{2}} u|^{2}dx+\tilde{\mu}+\sum_{j=1}^{\infty}\mu_{j}\delta_{x_{j}}$, for some $\tilde{\mu}\in {\mathcal{M}_{+}(\R^N)}$, $\{\mu_j\}\subseteq(0,1)$ and $\{x_j\}\subseteq\R^N$.
Up to subsequences, there exists a measure $\nu \in \mea$ such that $|u_{\eps}|^{2^{\ast}}\!dx \tows \nu$.
and  by Theorem \ref{the_cca} there exists a set of nonnegative numbers $\{\nu_j\}_{j\in J}$ such that (up to reordering the points $\{x_j\}$ and the $\{\mu_j\}$)
\begin{equation}\label{eq_cca}
\dys \nu=|u|^{2^\ast}\!dx+\sum_{j}\nu_{j}\delta_{x_{j}}  \ \ \text{and}\ \ \nu_{j}\leq S^{\ast}\mu_{j}^{\frac{2^{\ast}}{2}}.
\end{equation}

Using H\"older Inequality, we have 
\begin{eqnarray*}
\dys {F}_{\eps}(u_{\eps},\mu_{\eps}) \ = \ \int_{\Omega}|u_{\eps}|^{2^{\ast}-\eps}dx \ \leq \ \left(\int_{\Omega}|u_{\eps}|^{2^{\ast}}dx\right)^{\!\!\frac{2^{\ast}-\eps}{2^{\ast}}}\!|\Omega|^{\frac{\eps}{2^{\ast}}},
\end{eqnarray*}
hence, the definition of $\nu$ and~\eqref{eq_cca} yield
\begin{eqnarray*}
 \dys \limsup_{\eps\to 0^{+}}{F}_{\eps}(u_{\eps},\mu_{\eps})
& \leq &  \limsup_{\eps\to 0}\left(\int_{\Omega}|u_{\eps}|^{2^{\ast}}dx\right)^{\!\!\frac{2^{\ast}-\eps}{2^{\ast}}}\!|\Omega|^{\frac{\eps}{2^{\ast}}}\\
& \leq & \nu(\Omb) \ \leq \ \int_{\Omega}|u|^{2^{\ast}}dx + S^{\ast}\sum_{i=1}^{\infty}\mu_{i}^{\frac{2^{\ast}}{2}} \  \leq \ {F}(u,\mu).  
\end{eqnarray*}
\end{proof}

\vspace{2mm}

Now, we will prove the $\Gamma^{+}$-liminf inequality (ii).

It is convenient to define a relevant subset of configurations $\tilde{X}\subset X$ as follows
$$
\tilde{X}:=\Big\{ (u,\mu)\in \dot{H}^s(\Om)\!\times\!\mea : \mu=|(-\Delta)^{\frac{s}{2}}u|^2dx+\tilde{\mu}+\sum_{j=1}^n\mu_j\delta_{x_j}, \, \mu(\R^N)<1 \Big\}.
$$

For any pair $(u,\mu)$ in $\tilde{X}$ we will prove the existence of a recovery sequence $\{(\bu,\mbu)\}\subset X$ for the $\gamp$-liminf inequality, as stated in the following proposition.
\begin{proposition}\label{pro_recovery}
For any $(u,\mu)\in\tilde{X}$ there exists a sequence $\{(\bu,\mbu)\}\subset X$ such that $(\bu,\mbu)\towt(u,\mu)$ and
\begin{equation}\label{eq_recovery}
\dys \liminf_{\eps\to 0}F_\eps(\bu,\mbu)\, \geq\, F(u,\mu).
\end{equation}
\end{proposition}
\vspace{1mm}
Finally, we will prove the $\gamp$-liminf inequality in the whole space $X$, by a diagonal argument using recovery sequences for the elements of $\tilde{X}$.
\vspace{2mm}

In order to prove Proposition~\ref{pro_recovery}, first, for any point $x_j$ in $\Omb$ we construct a sequence $\{v^{j}_\eps\}$ that concentrates energy at $x_j$ (see forthcoming Proposition~\ref{pro_talenti}).
Then, we show that we can glue such sequences $\{v^{j}_\eps\}$ into a sequence $\{u^A_\eps\}$ such that it concentrates at any finite set of points $\{x_j\}$ in $\Omb$ (see Corollary~\ref{cor_talenti}).
Thus, the sequence $\{\ue^A\}$ will be the recovery sequence for a pair $(0,\mu) \in \tilde{X}$ when $\mu$ is purely atomic.

Finally, for any pair $(u,\mu)$ in $\tilde{X}$, we will able to join the function $u$ to the sequence $\{u^A_\eps\}$, adding suitably their corresponding measures, to obtain the desired recovery sequence~$\{(\bu,\mbu)\}$ satisfying~\eqref{eq_recovery}. More precisely, combining the lemmas in~Section~\ref{sec_lemma} with a careful choice of the supports of the approximating functions $\{u^A_\eps\}$ will give an admissible sequence~$\{(\bu,\mbu)\}\subseteq X$. Then, a precise calculation still based on Lemma~\ref{lem_cutoff} and Lemma~\ref{lem_commutator} will give~\eqref{eq_recovery}.

\vspace{1mm}

We start with the following result.

\begin{proposition}\label{pro_talenti}
For any $x_0\in \Omb$ there exists a sequence $\{v_\eps\}\subset \dot{H}^s(\Omega)$ such that
\begin{itemize}
\item[(i)]{
$\{(v_\eps, |(-\Delta)^{\frac{s}{2}}v_\eps|^2dx)\}\subseteq X$ and $\tau$-converges to $(0,\delta_{x_0})$ as $\eps\to 0$;\vspace{1mm}
}
\item[(ii)]{$\dys \lim_{\eps\to0}\text{\rm dist}_H\big(\text{\rm spt}\ v_\eps, \{x_0\}\big) = 0$;\vspace{1mm}
}
\item[(iii)]{$\dys \lim_{\eps\to0} \int_\Omega |v_\eps|^{2^{\ast}-\eps}dx = S^{\ast}$.
}
\end{itemize}
\end{proposition}
\vspace{1mm}
\begin{proof}
We assume that $x_0$ is an interior point of $\Om$ and we construct the sequence $\{v_\eps\}$ modifying the extremal functions $u$ for the Sobolev embedding $S^\ast$ given by Theorem~\ref{the_cotsiolis}.

Let $u\in \dot{H}^s(\R^N)$ defined as follows
$$
u(x)=\frac{c}{(1+|x-x_0|^2)^{\frac{N-2s}{2}}}, \ \ \forall x \in \R^N,
$$
where the positive constant $c$ is chosen such that $\|u\|_{\dot{H}^s}=1$.

If, for any positive $\eps$, we set $w_\eps(x):=\eps^{-\frac{N-2s}{2}}u(x/\eps)$, then we have
\begin{equation}\label{eq_s0}
\dys \int_{\R^N} |w_\eps|^{2^\ast}dx = S^\ast \ \ \text{and} \ \ \|w_\eps\|_{\dot{H}^s} = 1,
\end{equation}
by scaling invariance of $L^{2^{\ast}}\!$ and $\dot{H}^s$ norms.

Moreover, the function $w_\eps$ satisfies $\dys w_\eps \rightharpoonup 0$ in $L^{2^\ast}\!(\R^N)$ and $|(-\Delta)^{\frac{s}{2}}w_\eps|^2dx \tows \delta_{x_0}$ in $\mea$ as $\eps\to 0$, 
since a direct calculation for any $\rho>0$ gives
\begin{equation}\label{eq_h}
\dys
w_\eps \overset{\varepsilon \to 0}{\longrightarrow} 0 \ \text{in} \ L^{2^\ast}\!(\R^N\setminus \overline{B_\rho(x_0)}) \ \ \text{and} \ \ |(-\Delta)^{\frac{s}{2}}w_\eps|^2 \overset{\varepsilon \to 0}{\longrightarrow} 0 \ \text{in} \ L^2(\R^N\setminus  \overline{B_\rho(x_0)}).
\end{equation}

\vspace{1mm}

We want to localize the sequence $w_\eps$ in smaller and smaller neighborhoods of $x_0$.

For any fixed positive $\rho$, take a cut-off function $\varphi \in C^{\infty}_0(\R^N)$ such that $\varphi\equiv 1$ in $B_\rho(x_0)$, $\varphi\equiv 0$ in $\R^N\setminus B_{2\rho}(x_0)$ and $0\leq \varphi\leq 1$. For any $\eps>0$, we define
$$
\tilde{v}_\eps(x):=\varphi(x)w_\eps(x)
$$
and we claim that, as $\eps\to0$,
\begin{equation}\label{eq_convergenze}
\dys
\tilde{v}_\eps \rightharpoonup 0 \ \text{in} \ L^{2^\ast}\!(\Om), \ \ \
\|\tilde{v}_\eps\|_{\dot{H}^s} \to 1 \ \ \text{and} \ \
\dys
\int_\Om |\tilde{v}_\eps|^{2^\ast-\eps}dx \to S^\ast.
\end{equation}
The first convergence result in~\eqref{eq_convergenze} is a direct consequence of~\eqref{eq_h}. In order to prove the second, note that $\|u\|_{\dot{H}^s}=1$, hence Lemma~\ref{lem_cutoff} yields
\begin{eqnarray*}
\dys
\|\tilde{v}_\eps\|_{\dot{H}^s} & = & \int_{\R^N} |(-\Delta)^{\frac{s}{2}}\big(\varphi(x)w_\eps(x/\eps)\big)|^2dx \\
& = &
\int_{\R^N} |(-\Delta)^{\frac{s}{2}}\big(\varphi(\eps y) u (y)\big)|^2dy \  \overset{\varepsilon \to 0}{\longrightarrow} \ 1.
\end{eqnarray*}

The last convergence result in \eqref{eq_convergenze} is more delicate. We split the integral into two parts, namely $I_{1,\eps}$ and $I_{2,\eps}$ given by
$$
\dys
I_{1,\eps}:=\int_{\Om\cap\{w_\eps<1\}} |\varphi w_\eps|^{2^\ast-\eps}dx \ \ \ \text{and} \ \ \ I_{2,\eps}:=\int_{\Om\cap\{w_\eps\geq1\}}|\varphi w_\eps|^{2^\ast-\eps}dx.
$$
Since $|w_\eps|^{2^\ast-\eps}\leq 1$ in $\Om \cap \{w_\eps<1\}$ uniformly in $\eps$ and $|\varphi w_\eps(x)|^{2^\ast-\eps}\to 0$ a.e. as $\eps\to 0$, we deduce that $I_{1,\eps}$ vanishes as $\eps$ goes to 0.

For $I_{2,\eps}$  first we want to prove that
\begin{equation}\label{eq_s4}
\dys
\lim_{\eps\to 0} \| \frac{\varphi^{2^{\ast}-\eps}}{w_\eps^{\eps}}-1 \|_{L^{\infty}(\{w_{\eps} \geq 1\})} = 0.
\end{equation}
Note that, for $\eps$ small enough, we have $\{w_\eps\geq 1\} \subseteq B_\rho(x_0)$ and then 
$
\dys \frac{\varphi^{2^\ast-\eps}}{w_\eps^\eps}-1 = \frac{1}{w_\eps^\eps}-1$
 in $\Om\cap\{w_\eps<1\}$.
Hence, \eqref{eq_s4} follows once we show that
\begin{equation}\label{eq_s5}
\dys
\lim_{\eps\to 0} \| w_\eps^{\eps}-1 \|_{L^{\infty}(\{w_{\eps} \geq 1\})} = 0.
\end{equation}
Clearly, on $\Om\cap\{w_\eps\geq 1\}$ the function 
$w_\eps$ satisfies
$$
1\, \leq \, w_\eps^\eps \, \leq (\max w_\eps)^\eps \, = \, \left(c \eps^{-\frac{N-2s}{2}}\right)^\eps
$$ 
and thus we obtain~\eqref{eq_s5} and in turn~\eqref{eq_s4} as $\eps\to 0$.
\vspace{1mm}

Combining~\eqref{eq_s4} with~\eqref{eq_s0}, $I_{2,\eps}$ can be estimating as follows
\begin{eqnarray*}
\dys
I_{2,\eps} \! &= &\! \int_{\Om\cap\{w_\eps\geq 1\}} |\varphi w_\eps|^{2^\ast - \eps} dx 
\, = \, \int_{\Om\cap\{w_\eps\geq 1\}} |\frac{|\varphi|^{2^\ast-\eps}}{w_\eps^\eps}||w_\eps|^{2^\ast}dx \\
& = & \! \int_{\Om\cap\{w_\eps\geq 1\}} |w_\eps|^{2^\ast}dx + o(1) \,  \overset{\varepsilon \to 0}{\longrightarrow} \, S^\ast.
\end{eqnarray*}

\vspace{1mm}

Thus,~\eqref{eq_convergenze} holds for any $\rho>0$ small enough, whence a diagonal argument as $\rho \searrow 0$ gives a sequence $\{\tilde{v}_\eps\}$ such that~\eqref{eq_convergenze} holds, since $\tilde{v}_\eps\rightharpoonup 0$ in $L^{2^\ast}\!(\Om)$ and $\dys \lim_{\eps\to0}\textrm{dist}_H\big(\textrm{spt}\, \tilde{v}_\eps, \, \{x_0\}\big)=0.$

Note that, by Proposition~\ref{pro_supporto}, $\tilde{v}_\eps$ also satisfies
\begin{equation}\label{eq_j}
\dys
|(-\Delta)^{\frac{s}{2}}\tilde{v}_\eps|^2dx \to 0 \ \ \text{in} \ L^{2}(\R^N\setminus \overline{B_{\rho}(x_0)}) \ \ \text{as} \ \eps\to 0.
\end{equation}

\vspace{1mm}
Finally, for any $\eps>0$, we set 
\begin{equation*} 
v_\eps(x):= \frac{\tilde{v}_\eps(x)}{\|\tilde{v}_\eps\|_{\dot{H}^s}}.
\end{equation*}

Claim (i) follows readily  from~\eqref{eq_convergenze} and~\eqref{eq_j}. 
Claim~(ii) holds by construction, since the function $\tilde{v}_\eps$ has the same property. Finally, a simple calculation of the $L^{2^\ast\!-\eps}$ norm of the function $v_\eps$ gives
$$
\dys
\int_\Om |v_\eps|^{2^\ast-\eps}dx = \|\tilde{v}_\eps\|_{\dot{H}^s}^{-(2^\ast-\eps)}\!\int_\Om|\tilde{v}_\eps|^{2^\ast-\eps}dx \, \to \, S^\ast \ \ \text{as} \ \eps\to 0,
$$
which proves claim~(iii).
\vspace{1mm}

To complete the proof, we observe that the case of $x_0 \in \partial\Om$ can be obtained by a standard diagonal argument taking an approximating sequence of points $\{x_k\}\subseteq \Om$ converging to $x_0$ and the optimal sequences corresponding to each $x_k$.
\end{proof}

\begin{corollary}\label{cor_talenti}
For any finite set of distinct points $\{x_1, x_2, ..., x_n\}\subset \Omb$ and for any set of positive numbers $\{\mu_1, \mu_2, ..., \mu_n\}\subseteq\R$ such that $\sum_j\mu_j<1$, there exists a sequence $\{u_\eps^A\}\subset \dot{H}^s(\Om)$ such that
\begin{itemize}
\item[(i)]{
$\{(u^A_\eps, |(-\Delta)^{\frac{s}{2}}u_\eps^A|^2dx)\}\subseteq X$ and $\tau$-converges to $(0, \sum_{j=1}^{n}\mu_j\delta_{x_j})$ as $\eps\to 0$;\vspace{1mm}
}
\item[(ii)]{${\dys \lim_{\eps\to0}}\,\text{\rm dist}_H\Big(\text{\rm spt}\ u^A_\eps\,, \ \bigcup_{j}\{x_j\}\Big) = 0$.\vspace{1mm}
}
\item[(iii)]{$\dys \lim_{\eps\to0} \int_\Omega |u^A_\eps|^{2^{\ast}-\eps}dx = S^{\ast}\sum_{j=1}^{n}\mu_j^{\frac{2^{\ast}}{2}}$.
}
\end{itemize}
\end{corollary}
\vspace{1mm}
\begin{proof}
Let us set $A_{j}:=B_{r_{j}}(x_{j})\cap\Om$ for any $j=1, 2 ..., n$, with radii  
$r_{j}$ and $r_{i}$ such that $\text{dist}(A_{j},A_{i})>0$.
By Proposition~\ref{pro_talenti}, there exists a sequence $\{(v_\eps^{j},\mu_\eps^{j})\}\subset X$, with $\mu_\eps^{j}=|(-\Delta)^{\frac{s}{2}}v_\eps^{j}|^{2}dx$ such that $(\ue^{j},\mue^{j}) \towt (0,\delta_{x_{i}})$, 
spt$\, v_\eps^j \subset A_j$, 
dist$_H\big(\text{spt} \ v_\eps^j,\{x_j\}\big)\to 0$ as $\eps\to0$ and
\begin{equation*} 
\dys \lim_{\eps\to0} \int_\Omega |v_\eps^j|^{2^{\ast}-\eps}dx = S^{\ast},  \ \text{for} \ j=1, 2, ..., n. 
\end{equation*}

Let us set $\dys \ue^A:=\sum_{j=1}^{n}\sqrt{\mu_{j}} v_\eps^{j}$.

Estimating the energy of the sequence $\{|(-\Delta)^{\frac{s}{2}}\ue^A)|^2dx\}$ gives
\begin{eqnarray}\label{eq_enmue}
\dys \int_{\R^N}|(-\Delta)^{\frac{s}{2}}\ue^A|^2dx & = & \sum_{j=1}^n\mu_j\!\int_{\R^N}|(-\Delta)^{\frac{s}{2}}v^j_\eps|^2dx \nonumber \\
& & + 2\sum_{i, j=1, i<j}^n\sqrt{\mu_i\mu_j}\langle(-\Delta)^{\frac{s}{2}}v_\eps^i, (-\Delta)^{\frac{s}{2}}v^j_\eps\rangle_{L^2(\R^N)}.
\end{eqnarray}
We claim that the last sum  in the  formula above converges to zero as $\eps$ goes to zero. Indeed, by Cauchy-Schwarz inequality, we get
\begin{eqnarray}\label{eq_doppiop}
&& \dys 
\left|\langle(-\Delta)^{\frac{s}{2}}v^i_{\eps}, (-\Delta)^{\frac{s}{2}}v^j_{\eps}\rangle_{L^2(\R^N)}\right| \nonumber \\
&& \qquad\qquad\qquad \leq  \left(\int_{\R^N}|(-\Delta)^{\frac{s}{2}}v^i_{\eps}|^2dx\right)^{\!\frac{1}{2}}\left(\int_{H_i}|(-\Delta)^{\frac{s}{2}}v^j_{\eps}|^2dx\right)^{\!\frac{1}{2}} \nonumber \\
& & \qquad\qquad\qquad  \quad \ +\left(\int_{\R^N}|(-\Delta)^{\frac{s}{2}}v^j_{\eps}|^2dx\right)^{\!\frac{1}{2}}\left(\int_{H_j}|(-\Delta)^{\frac{s}{2}}v^i_{\eps}|^2dx\right)^{\!\frac{1}{2}}\!,
\end{eqnarray}
where, for $i$ and $j$ fixed, we have divided the whole space $\R^N$ into two complementary half-spaces $H_i$ and $H_j$ such that $A_i\subset H_i$ and $A_j\subset H_j$.

Note that $\int|(-\Delta)^{\frac{s}{2}}v^j_{\eps}|^2dx$ is smaller than 1 uniformly with respect to $\eps$ because $\{(v_\eps^j, |(-\Delta)^{\frac{s}{2}} v_\eps^j|^2dx)\} \subseteq X$. Thus, (\ref{eq_doppiop}) becomes
\begin{eqnarray*}
&& \dys 
\left|\langle(-\Delta)^{\frac{s}{2}}v^i_{\eps}, (-\Delta)^{\frac{s}{2}}v^j_{\eps}\rangle_{L^2(\R^N)}\right|  \\
&& \qquad\qquad\qquad\qquad\leq   \left(\int_{H_i}|(-\Delta)^{\frac{s}{2}}v^j_{\eps}|^2dx\right)^{\!\frac{1}{2}} +\left(\int_{H_j}|(-\Delta)^{\frac{s}{2}}v^i_{\eps}|^2dx\right)^{\!\frac{1}{2}}\!\!.
\end{eqnarray*}

On the other hand, since the measure $|(-\Delta)^{\frac{s}{2}}v^j_{\eps}|^2dx$ converges to $\delta_{x_j}$ in~$\mea$ as $\eps\to 0$ and spt$\, v_\eps^j \subseteq A_j$ for all $j=1, 2, ..., n$, Proposition \ref{pro_supporto} yields
$$
\dys \int_{H_i}|(-\Delta)^{\frac{s}{2}}v^j_{\eps}|^2dx \to 0 \ \ \text{as}\  \eps\to 0 \ \text{for} \ i\neq j,
$$
that in turn implies
\begin{equation}\label{eq_supp}
\dys \langle(-\Delta)^{\frac{s}{2}}v^i_{\eps}, (-\Delta)^{\frac{s}{2}}v^j_{\eps}\rangle_{L^2(\R^N)} \to 0 \ \ \text{as} \ \eps\to 0,
\end{equation}
Combining (\ref{eq_enmue}) and (\ref{eq_supp}) with the fact that each $v^j_\eps$ concentrates energy at $x_j$, in the sense of Proposition~\ref{pro_talenti}, we deduce that the constructed sequence $\{|(-\Delta)^{\frac{s}{2}}\ue^A|^2dx\}$ satisfies
$$
|(-\Delta)^{\frac{s}{2}}\ue^A|^2dx \tows \sum_{j=1}^n\mu_j\delta_{x_j} \ \, \text{in} \ \mea.
$$

Finally, since $\sum_j\mu_j<1$, by~\eqref{eq_enmue} we also deduce that
$\int_{\R^N}|(-\Delta)^{\frac{s}{2}}\ue^A|^2dx\leq 1$, for $\eps$ small, hence 
$\{(\ue^A,|(-\Delta)^{\frac{s}{2}}\ue^A|^2dx)\}\subset X$, for $\eps$ small and claim~(i) is completely proved.

\vspace{1mm}

Note that (ii) follows by construction, because of Proposition~\ref{pro_talenti}.

\vspace{1mm}

Moreover, 
since spt$\, v_\eps^j$ are mutually disjoint and $v^j_\eps$ satisfies Proposition~\ref{pro_talenti} (iii), 
we have
$$
\int_{\Om}|\ue^A|^{2^{\ast}\!-\eps}dx\,=\,\sum_{j=1}^n\mu_{j}^{\frac{2^{\ast}\!-\eps}{2}}\int_{A_{j}}|v_\eps^{j}|^{2^{\ast}\!-\eps}dx \ \overset{\eps\to 0}{\longrightarrow} \ S^{\ast}\sum_{j=1}^n\mu_{j}^{\frac{2^{\ast}}{2}},
$$
which concludes the proof of claim (iii). 
\end{proof}

\vspace{1mm}

Now, we are in position to prove the $\gamp$-liminf inequality for the set of configurations $\tilde{X}$ as stated in Proposition~\ref{pro_recovery}. The main contribution is given by the sequence $\{(\ue^A,|(-\Delta)^{\frac{s}{2}}\ue^A|^2dx)\}$ built in Corollary~\ref{cor_talenti}, but we have to carefully modify it in order to obtain the desired recovery sequence $\{(\bu, \mbu)\}$.

\begin{proof}[\bf Proof of Proposition~\ref{pro_recovery}]
Let $(u,\mu)$ be any fixed pair in $\tilde{X}$, i.~\!e., $u\in \dot{H}^s(\Om)$ and $\dys \mu=|(-\Delta)^{\frac{s}2}u|^2dx+\tilde{\mu}+\sum_{j=1}^n\mu_j\delta_{x_j} \in \mea$,
with $\mu(\R^N)<1$ and let
 $\{\ue^A\}$ be the corresponding sequence given by Corollary~\ref{cor_talenti}.
 
 \vspace{2mm}

For $\sigma>0$, take a cut-off function $\varphi_{\sigma}$ in $C^{\infty}_0(\R^N)$ such that $\varphi_\sigma\equiv 0$ in $B_{\rho_\sigma}(x_j)$, for $j=1,2,...,n$, $\varphi_\sigma\equiv 1$ in $\Om\setminus\bigcup_j B_{2\rho_\sigma}(x_j)$, with $\rho_\sigma\to 0$ as $\sigma\to 0$, $\dys \varphi_\sigma = 1 - \sum_{j=1}^n\bar{\varphi}\big(\frac{x-x_j}{\rho_\sigma}\big)$, $\bar{\varphi} \in C^\infty_0(B_2)$, $\bar{\varphi}\equiv 1$ on $\overline{B}_1$, $0\leq \overline{\varphi}\leq 1$.
\vspace{2mm}

Now, we can define the sequence $\{(\bu, \mbu)\}$ as follows
\[ \dys \bu \!\! =  \!\! \bar{u}_{\eps,\sigma}:=u\varphi_\sigma +\ue^A \quad \, , \quad \, \dys \mbu \!\!  =  \!\!\bar{\mu}_{\eps,\sigma}:=\tilde{\mu}+|(-\Delta)^{\frac{s}{2}}(u\varphi_\sigma+\ue^A)|^2dx \]
and we claim that this is a recovery sequence for $(u\varphi_\sigma, |(-\Delta)^{\frac{s}{2}}(u\varphi_\sigma)|^2dx+\tilde{\mu}+\sum_{j=1}^n\mu_j\delta_{x_j})$ as $\eps \to 0$.
\label{rem_diag}
Note that we will play with two positive parameters, namely $\eps$ (which is the parameter for the atomic part of $\mu$) and $\sigma$ (which will control the diffuse part of $\mu$). We will take limits in these parameters in the following order: first $\eps\to 0$, then $\sigma\to 0$. The recovery sequence for $(u,\mu)$ will be actually given by a further diagonal argument.
\vspace{1mm}

First, we claim that $\{(\bu, \mbu)\}\subset X$ for $\eps$ and $\sigma$ small enough. Since we have $\bu\in \dot{H}^s(\Om)$ (because $\varphi_\sigma$ is a multiplier in $\dot{H}^s(\Om)$; see~\cite{MS}) and $\mbu \geq |(-\Delta)^{\frac{s}{2}}\bar{u}_\eps|^2dx$, this claim reduces to proving that
\begin{equation}\label{eq_4d0}
\mbu(\R^N)\leq 1.
\end{equation}
In order to check~\eqref{eq_4d0}, for any $\eps, \sigma >0$ we compute
\begin{eqnarray}\label{eq_4d}
\dys
\mbu(\R^N) \! & = & \! \tilde{\mu}(\R^N) + \int_{\R^N} |(-\Delta)^{\frac{s}{2}}(u\varphi_\sigma+\ue^A)|^2dx \nonumber \\
& = & \! \tilde{\mu}(\R^N) + \int_{\R^N} |(-\Delta)^{\frac{s}{2}}(u\varphi_\sigma)|^2dx 
+ \int_{\R^N} |(-\Delta)^{\frac{s}{2}}(\ue^A)|^2dx \\
& & \! +\, \langle(-\Delta)^{\frac{s}{2}}(u\varphi_\sigma), \, (-\Delta)^{\frac{s}{2}}\ue^A\rangle_{L^2(\R^N)}. \nonumber
\end{eqnarray}
We can treat the last three terms in the right-handside of equation~\eqref{eq_4d} as follows.

For $\sigma>0$ fixed, Corollary~\ref{cor_talenti}-(i) and Proposition~\ref{pro_supporto} yield
\begin{equation}\label{eq_4d2}
\dys
\lim_{\eps\to 0}\int_{\R^N}|(-\Delta)^{\frac{s}{2}}\ue^A|^2dx \, = \, \sum_{j=1}^n \mu_i.
\end{equation}

Again by Corollary~\ref{cor_talenti} $\ue^A \rightharpoonup 0$ in $\dot{H}^s(\Om)$, hence we have
\begin{equation}\label{eq_4d3}
\dys
\lim_{\eps\to 0}\langle(-\Delta)^{\frac{s}{2}}(u\varphi_\sigma), \, (-\Delta)^{\frac{s}{2}}\ue^A\rangle_{L^2(\R^N)} \, = \, 0.
\end{equation}

Finally, from the definition of $\varphi_\sigma$ and Lemma~\ref{lem_cutoff}, we have
\begin{equation}\label{eq_4dlast}
\dys
\lim_{\sigma\to 0}\int_{\R^N}|(-\Delta)^{\frac{s}{2}}(u\varphi_\sigma)|^2dx \, = \, \int_{\R^N} |(-\Delta)^{\frac{s}{2}}u|^2dx.
\end{equation}

Thus, combining~\eqref{eq_4d2}, \eqref{eq_4d3}, \eqref{eq_4dlast} with the fact that $\mu(\R^N)$ is strictly less than 1 (recall that $(u,\mu)\in \tilde{X}$), we can deduce the inequality in~\eqref{eq_4d0} for $\eps$ and $\sigma$ small enough.

\vspace{2mm}

We prove that $\{(\bu\mbu)\}$ $\tau$-converges to $(u,\mu)$, i.e., 
\begin{equation}\label{eq_5d}
\dys \bu\rightharpoonup u \ \text{in} \ L^{2^\ast}\!(\Om) \ \ \ \text{and} \ \ \ \dys \mbu\tows \mu \ \text{in} \ \mea.
\end{equation}

Clearly, $|u\varphi_\sigma - u|^{2^\ast} \, = \,  |1-\varphi_\sigma|^{2^\ast}|u|^{2^\ast} \, \leq \, |u|^{2^\ast}$,
thus, $\{u\varphi_\sigma\}$ converges strongly to $u$ in $L^{2^\ast}$, as $\sigma\to 0$, by Lebesgue's Dominated Convergence Theorem, and then the first convergence result in~\eqref{eq_5d} follows from the fact that the sequence $\{\ue^A\}$ weakly converges to $0$ in $L^{2^\ast}\!(\Om)$ as $\eps$ goes to $0$.
\vspace{2mm}

The second convergence result in~\eqref{eq_5d} is a consequence of the convergence in the sense of measures of the sequence $\{|(-\Delta)^{\frac{s}{2}}\ue^A|^2dx\}$ to the finite sum of Dirac masses~$\sum_j\mu_j\delta_{x_j}$,
together with the fact that $\ue^A\rightharpoonup 0$ as $\eps\to 0$ and $u\varphi_\sigma \to u$ in $\dot{H}^s(\Om)$ as $\sigma\to 0$ by Corollary~\ref{cor_talenti} and Lemma~\ref{lem_cutoff}, respectively. 
Indeed, by arguing as in~\eqref{eq_4d}, \eqref{eq_4dlast} and \eqref{eq_4d2}, for any $\psi\in C^0_0(\R^N)$, we have
\begin{eqnarray*}
\dys
\lim_{\sigma\to0} \lim_{\eps\to0} \int_{\R^N}\psi d\mbu  & = &  \lim_{\sigma\to0} \lim_{\eps\to0}\int_{\R^N}\psi d\tilde{\mu} + \int_{\R^N}\psi|(-\Delta)^{\frac{s}{2}}(u\varphi_\sigma+\ue^A)|^2dx \\
& = &  \int_{\R^N}\psi d\tilde{\mu} + \int_{\R^N}\psi|(-\Delta)^{\frac{s}{2}}u|^2dx + \int_{\R^N}\psi d\sum_j^n\mu_j\delta_{x_j} \\
& = &  \int_{\R^N}\psi d\mu,
\end{eqnarray*}
that completely proves~\eqref{eq_5d}.
\vspace{2mm}

In order to complete the proof, it remains to show that its energy $F_\eps(\bu,\mbu)$ satisfies the liminf inequality stated in~\eqref{eq_recovery}. Since $\text{dist}\Big(\text{spt}\, (u\varphi_\sigma), \, \bigcup_j B_{\rho_\sigma}(x_j) \Big)>0$, we can split the integral in $F_\eps(\bu,\mbu)$ as follows
\begin{eqnarray}\label{eq_7d}
\dys
\\
F_\eps(\bu,\mbu) \ = \ \int_\Om |u\varphi_\sigma + \ue^A|^{2^{\ast}-\eps}dx \ = \ \int_\Om |u\varphi_\sigma|^{2^{\ast}-\eps}dx + \int_\Om |\ue^A|^{2^{\ast}-\eps}dx.\nonumber
\end{eqnarray}
By Dominated Convergence Theorem, we have
\begin{equation}\label{eq_8d}
\lim_{\sigma\to0} \lim_{\eps\to0} \int_\Om |u\varphi_\sigma|^{2^{\ast}-\eps}dx \, = \,
\int_\Om |u|^{2^\ast}dx.
\end{equation}
On the other hand, taking Corollary~\ref{cor_talenti}-(iii) into account, we have
\begin{equation}\label{eq_8d2}
\dys
\int_\Om|\ue^A|^{2^{\ast}-\eps} \, \to \, S^{\ast}\sum_{j=1}^n\mu_j^{\frac{2^{\ast}}{2}} \ \, \text{as} \ \eps\to 0.
\end{equation}
Finally, combining~\eqref{eq_7d}, \eqref{eq_8d} and \eqref{eq_8d2}, we obtain (up to the diagonal argument on $\eps$ and $\sigma$ mentioned at page~\pageref{rem_diag})
$$
\dys
\qquad\qquad \liminf_{\eps \to 0} F_\eps(\bu,\mbu) \, = \, \int_\Om|u|^{2^\ast}dx + S^\ast \sum_{j=1}^n \mu_i^{\frac{2^\ast}{2}}
\, = \, F(u,\mu).  
$$
\end{proof}

\begin{proof}[\bf Completion of the proof of the $\gamp$-liminf inequality]
In view of Proposition~\ref{pro_recovery}, the $\gamp$-liminf inequality in Theorem~\ref{the_gamma-intro} holds for any $(u,\mu) \in \tilde{X}$. Thus, it is enough to check that $\tilde{X}\subseteq X$ is $\tau$-sequentially dense by an explicit approximation and that $F$ is continuous with respect to this approximation, in order to conclude by a standard diagonal argument.
\vspace{1mm}

For any pair $(u,\mu)\in X$, we consider the sequence $\{(u_n, \mu_n)\}$ defined as
$$
\dys u_n :=  c_n u \ \ \ \text{and} \ \ \
\dys \mu_n  :=  c^2_n|(-\Delta)^{\frac{s}{2}}u|^2dx + c^2_n\tilde{\mu} +c^2_n\sum_{j=1}^n \mu_j\delta_{x_j},
$$
where $\{c_n\}\subset (0,1)$ is any increasing sequence such that $c_n \nearrow 1$ as $n\to \infty$.\vspace{1mm}

Clearly, the sequence $\{(u_n,\mu_n)\}$ is in $\tilde{X}$, since, for any $n\in \N$, $u_n \in \dot{H}^s(\Om)$ and $\mu_n$ is a measure with a finite number of atoms such that $\dys \mu_n(\R^N) \, \leq \, c^2_n\,\mu(\R^N) \, \leq \, c^2_n \, < 1$. Moreover, $(u_n, \,u_n) \towt (u,\mu)$ as $n\to \infty$, because $u_n \to u$ in $\dot{H}^s(\Om)$ (hence weakly in $L^{2^\ast}\!(\Om)$) and, for any $\psi \in C^0_0(\R^N)$,
\begin{eqnarray*}
\dys
&& \hspace{-1cm} \int_{\R^N} \psi d\mu_n \, = \, \int_{\R^N} \psi c^2_n|(-\Delta)^{\frac{s}{2}}u|^2dx + \int_{\R^N}\psi c^2_n d\tilde{\mu} + c^2_n \sum_{j=1}^n \mu_j\psi(x_j) \\
&& \qquad = \, c^2_n\left(\int_{\R^N} \psi |(-\Delta)^{\frac{s}{2}}u|^2dx + \int_{\R^N}\!\psi  d\tilde{\mu}+ \sum_{j=1}^n \mu_j\psi(x_j)\right) \, \overset{n\to \infty}{\longrightarrow} \, \int_{\R^N}\!\psi d\mu.
\end{eqnarray*}
\vspace{1mm}

Finally, evaluating the  functional $F$, we have
\begin{eqnarray*}
\dys
F(u_n, \mu_n) & = & \int_\Om c^{2^\ast}_n |u|^{2^\ast}dx + S^{\ast}\sum_{j=1}^n(c^2_n\mu_j)^{\frac{2^\ast}{2}} \\
& = &  c^{2^\ast}_n \left(\int_\Om|u|^{2^\ast}dx + S^{\ast}\sum_{j=1}^n \mu_j^{\frac{2^\ast}{2}}\right) \ \to \ F(u,\mu), \ \ \text{as} \ n \to \infty.  
\end{eqnarray*}
\end{proof}

\subsection{Concentration of optimal functions} 
Here we show that, due to the $\gamp$-convergence result,  the   maximizers $\{\ue\}$ for the variational problem (\ref{problema}) concentrate energy at one point $x_{0} \in \Omb$ when $\eps$ goes to zero.  The key result is the following optimal upper bound for the limit functional $F$ on the space X.
\begin{lemma}\label{lem_sstar}
For every $(u,\mu) \in X$, we have
\begin{equation}\label{eq_soboc}
{F}(u,\mu) \leq S^{\ast}
\end{equation}
and the equality holds if and only if $(u,\mu)=(0,\delta_{x_0})$ for some $x_0 \in \Omb$.
\end{lemma}
\vspace{1mm}
\begin{proof}
We adapt the argument in the proof of~\cite[Lemma~3.6]{amargarroni} for the case $s=1$, using the well-known convexity trick by P.~\!L. Lions. 
\vspace{1mm}

For every $(u,\mu)\in X$, by Sobolev inequality (\ref{eq_sobolev1}), we have
$$
{F}(u,\mu)  \equiv  \int_{\Om}|u|^{2^{\ast}}\!dx+S^{\ast}\sum_{j=1}^{\infty}\mu_{j}^{\frac{2^{\ast}}{2}} \leq  S^{\ast}\left(\int_{\R^N}|(-\Delta)^{\frac{s}{2}} u|^{2}dx\right)^{\!\!\frac{2^{\ast}}{2}}+S^{\ast}\sum_{j=1}^{\infty}\mu_{j}^{\frac{2^{\ast}}{2}}.
$$
Now, by the convexity of the function $t\mapsto t^{\frac{2^{\ast}}{2}}$, for every fixed $s \in (0,N/2)$, we get
\begin{eqnarray}\label{numero}
{F}(u,\mu)  & \leq & S^{\ast}\left(\int_{\R^N}|(-\Delta)^{\frac{s}{2}} u|^{2}dx\right)^{\!\!\frac{2^{\ast}}{2}}+S^{\ast}\sum_{j=1}^{\infty}\mu_{j}^{\frac{2^{\ast}}{2}} \nonumber \\ 
& \leq &  S^{\ast}\left(\int_{\R^N}|(-\Delta)^{\frac{s}{2}} u|^{2}dx+\sum_{j=1}^{\infty}\mu_{j}\right)^{\!\!\frac{2^{\ast}}{2}}  \leq \, S^{\ast}\!\left(\mu(\R^N)\right)^{\!\frac{2^\ast}{2}} \, \leq \, S^{\ast}  
\end{eqnarray}
which proves~\eqref{eq_soboc}.

\vspace{1mm}
Note that equality clearly holds if $(u,\mu)=(0,\delta_{x_{0}}),$ for some $x_{0} \in \Omb$.

Assume that equality in~\eqref{eq_soboc} holds for some pair $(u,\mu) \in X$. Then, each inequality in~\eqref{numero} is in effect an equality. In particular, we deduce $\tilde{\mu}=0$. If $u\neq 0$ then we also deduce by convexity that $\mu_j=0$ for every $j$. In turn, this fact yields $\mu=|(-\Delta)^{\frac{s}{2}}u|^2dx$ and $u\in \dot{H}^s(\Om)$ is optimal in Sobolev inequality~\eqref{eq_sobolev1}, which contradicts Theorem~\ref{the_cotsiolis}. Thus, $u=0$, equation~\eqref{numero} and the strict convexity implies that $\mu=\delta_{x_0}$ for some $x_0\in \Omb$ as claimed. 
\end{proof}
\vspace{1mm}

\begin{proof}[\bf Proof of Corollary~\ref{cor_concentration}]
Now, by Theorem \ref{the_gamma-intro} and standard  $\gamp$-convergence properties, it follows that every sequence of maximizers of $F_\eps$, which is in the form  $\{(\ue, |(-\Delta)^{\frac{s}{2}}\ue|^2dx)\}$ in view of Proposition~\ref{pro_argmax}, must converge (up to subsequences) to a pair $(u,\mu)\in X$ which is a maximizer for $F$, i.e.
$$
\dys (\ue,|(-\Delta)^{\frac{s}{2}}\ue|^2dx)\towt(u,\mu), \ \ \ \text{with} \ F(u,\mu)=\max_{X(\overline{\Om})}F.
$$
By Lemma \ref{lem_sstar}, we have the upper bound $\dys F(u,\mu)\leq S^{\ast}$  for every $(u,\mu)\in X$ 
and the equality is achieved if and only if $(u,\mu)=(0,\delta_{x_0})$ for some $x_0\in \Omb$. Hence, it follows that
$\dys
(\ue, |(-\Delta)^{\frac{s}{2}}\ue|^2dx)\towt(0,\delta_{x_0})$, which is the desired concentration property for the energy density.
\end{proof}

\vspace{2mm}

\end{document}